\documentclass[12pt]{article}
\usepackage{amsfonts,amssymb,amsmath,eucal}
\usepackage{oldgerm}
\usepackage{eufrak}
\usepackage[dvips]{graphics}

\newcommand{\ffP}{\mathfrak{P}}
\newcommand{\gl}{\mathfrak{gl}}

\newcommand{\Z}{\mathbb{Z}}
\newcommand{\bb}{\mathsf{b}}
\newcommand{\zz}{\mathfrak{z}}
\newcommand{\bT}{\mathbb{T}}
\newcommand{\ttC}{\mathtt{C}}
\newcommand{\fb}{\mathbf{f}}
\newcommand{\gb}{\mathbf{g}}
\newcommand{\pb}{\mathbf{p}}
\newcommand{\pbb}{\mathbf{\bar p}}
\newcommand{\eb}{\mathbf{e}}
\newcommand{\C}{\mathbb{C}}
\newcommand{\R}{\mathbb{R}}
\newcommand{\Q}{\mathbb{Q}}
\newcommand{\sQ}{\mathsf{Q}}
\newcommand{\QM}{\mathcal{QM}}
\newcommand{\cM}{\mathcal{M}}
\newcommand{\cE}{\mathcal{E}}
\newcommand{\cP}{\mathcal{P}}
\newcommand{\cPb}{\mathcal{\bar P}}
\newcommand{\cC}{\mathcal{C}}
\newcommand{\cG}{\mathcal{G}}
\newcommand{\bW}{\mathsf{w}}
\newcommand{\oW}{\textfrak{W}}
\newcommand{\Lb}{\bar\Lambda}
\newcommand{\vac}{v_\emptyset}
\newcommand{\LV}{{\Lambda^{\frac{\infty}2}}_0V}
\renewcommand{\th}{\vartheta}

\newcommand{\rang}{\right\rangle}
\newcommand{\lang}{\left\langle} 
\newcommand{\ve}[1]{\left|#1\right\rangle}

\DeclareMathOperator{\Aut}{Aut}
\DeclareMathOperator{\wt}{wt}
\DeclareMathOperator{\tr}{tr}

\newtheorem{Theorem}{Theorem}
\newtheorem{Lemma}{Lemma}

\newtheorem{Proposition}[Lemma]{Proposition}

\begin{document}

\title{Pillowcases and quasimodular forms}

\author{Alex Eskin\footnote{Department of Mathematics,
University of Chicago, 5734 South University Avenue
Chicago, Illinois 60637} \,\, and 
\setcounter{footnote}{6}
Andrei Okounkov\footnote{Department of Mathematics, Princeton University, 
Fine Hall, Washington Road, Princeton  New Jersey 08544\newline 
The authors were partially supported by NSF and Packard foundation}}

\date{\emph{To Vladimir Drinfeld on his 50th birthday}}

\maketitle 

\begin{abstract}
 We prove that natural generating functions for 
enumeration of branched coverings of the pillowcase orbifold
are level $2$ quasimodular forms. This gives a way to 
compute the volumes of the strata of the moduli 
space of quadratic differentials. 
\end{abstract}

\section{Introduction}

\subsection{Pillowcase covers and quadratic differentials}

Consider a complex torus $\bT^2=\C/L$, where $L\subset \C$
is a lattice. Its quotient 
$$
\ffP = \bT^2/\pm 
$$
by the automorphism $z\mapsto -z$ is a sphere with four 
$(\Z/2)$-orbifold points which is sometimes 
called the \emph{pillowcase}
orbifold. The map $\bT^2\to\ffP$ is essentially the 
Weierstra\ss\ $\wp$-function. 
The quadratic differential $(dz)^2$ on $\bT^2$ 
descends to a quadratic differential on $\ffP$. 
Viewed as a quadratic differential on the Riemann sphere, 
$(dz)^2$ has simple poles at corner points. 

Let $\mu$ be a partition and $\nu$ a partition of an even 
number into \emph{odd} parts. We are interested in 
enumeration of degree $2d$ maps 
\begin{equation}
  \label{cover_pi}
   \pi: \cC\to \ffP
\end{equation}
with the following ramification data. Viewed as a map
to the sphere, $\pi$ has profile $(\nu,2^{d-|\nu|/2})$ over 
$0\in\ffP$ and profile $(2^{d})$ over the other three 
corners of $\ffP$. Additionally, $\pi$ has profile 
$(\mu_i,1^{2d-\mu_i})$ over some $\ell(\mu)$ given points of 
$\ffP$ and unramified elsewhere. Here $\ell(\mu)$ is the 
number of parts in $\mu$. This ramification data
determines the genus of $\cC$ by
$$
\chi(\cC) = \ell(\mu)+\ell(\nu) - |\mu|-|\nu|/2\,. 
$$
In principle, one could allow more general ramifications
over $0$ and the nonorbifold points, but this more general 
problem is readily reduced to the one above\footnote{ 
{}From first principles, the count of the branched coverings
does not change if one replaces two ramification conditions
by the product of the corresponding conjugacy classes in 
the class algebra of the symmetric group. In this way, one can 
generate complicated ramifications from simpler ones.}.

Pulling back $(dz)^2$ via $\pi$ gives a quadratic 
differential on $\cC$ with zeros of multiplicities
$\{\nu_i-2\}$ and $\{2\mu_i-2\}$. The periods of this 
differential, by construction, lie in a translate of 
a certain lattice. The enumeration of covers
$\pi$ is thus related to lattice point enumeration 
in the natural strata of the
\emph{moduli space of quadratic differentials}.  
In particular, the $d\to\infty$
asymptotics gives the volumes of these strata. 
These volumes are of considerable interest in ergodic
theory, in particular in connection with billiards in 
rational polygons, see \cite{EMZ,MZ}. Their computation was
the main motivation for the present work. 

A different way to compute the volume of 
the principal stratum was found by M.~Mirzakhani \cite{Mir}. 

We thank H.~Masur, M.~Mirzakhani, N.~Nekrasov, R.~Pandharipande, and 
A.~Zorich for discussions and the anonymous referee for helpful 
suggestions.

\subsection{Generating functions}

\subsubsection{}

Two covers $\pi_i: \cC_i\to \ffP$, $i=1,2$, are identified if there
is an isomorphism $f: C_1 \to C_2$ such that $\pi_1=f\circ \pi_2$. 
In particular, associated to every cover $\pi$ is a finite
group $\Aut(\pi)$. This group is trivial for most connected
covers, see e.g.\
Section 3.1 in \cite{EO}. We form the following generating function 
\begin{equation}
  \label{Z_mu_nu}
   Z(\mu,\nu;q) = \sum_\pi \frac{q^{\deg \pi}}{|\Aut(\pi)|} \,,
\end{equation}
where $\pi$ ranges over all inequivalent covers \eqref{cover_pi}
with ramification data $\mu$ and $\nu$ as above. Note that 
the degree of any such $\pi$ is even. 

In particular, for $\mu=\nu=\emptyset$ any connected cover
has the form 
$$
\pi: \bT^2 \xrightarrow{\pi'} \bT^2 \to \bT^2/\pm 
$$
with $\pi'$ unramified. We have $|\Aut(\pi)|=2|\Aut(\pi')|$
corresponding to the lift of $\pm$. From the enumeration of
possible $\pi'$ we obtain,  
$$
Z(\emptyset,\emptyset;q) = \prod_{n} (1-q^{2n})^{-1/2} \,.
$$
By definition, we set
\begin{equation}
  \label{Z'}
    Z'(\mu,\nu;q) = \frac{Z(\mu,\nu;q)}{Z(\emptyset,\emptyset;q)} \,.
\end{equation}
This enumerates covers without unramified connected components. 
By the usual inclusion-exclusion, one can extract from \eqref{Z'}
a generating function for connected covers. This generating 
function for connected covers will be denoted by $Z^\circ(\mu,\nu;q)$.

\subsubsection{}

Recall the classical level 1 Eisenstein series
$$
E_{2k}(q) = \frac{\zeta(1-2k)}{2} + \sum_{n=1}^\infty \left(
\sum_{d|n} d^{2k-1}\right) q^n \,, \quad k=1,2,\dots \,.
$$
The algebra they generate is called the algebra $\QM(\Gamma(1))$
of \emph{quasimodular forms} for $\Gamma(1)=SL_2(\Z)$, see \cite{KZ}
and also below in Section \ref{squasi}. It is known 
that $E_2$, $E_4$, and $E_6$ are free 
commutative generators of $\QM(\Gamma(1))$. The algebra $\QM(\Gamma(1))$
is naturally graded by weight, where $\wt E_{2k} = 2k$. 
Clearly, for any integer $N$, $E_{2k}(q^N)$ is a quasimodular form 
of weight $2k$ 
for the group 
$$
\Gamma_0(N)= \left\{
  \begin{pmatrix}
    a & b \\
  c & d
  \end{pmatrix}\,\Big|\,\, c\equiv 0 \bmod N 
\right\} \subset SL_2(\Z) \,.
$$

The quasimodular forms that will appear in this paper
will typically be inhomogeneous, so instead of weight grading 
we will only keep track of the corresponding 
filtration. We define the weight of a partition by
$$
\wt \mu = |\mu| + \ell(\mu)
$$
The main result of this paper is the following 

\begin{Theorem}\label{T1} The series $Z'(\mu,\nu;q)$
is a polynomial in $E_2(q^2)$, $E_2(q^4)$, and $E_4(q^4)$
of weight $\wt \mu + |\nu|/2$. 
\end{Theorem}

Several explicit examples of the forms $Z'(\mu,\nu;q)$ 
are given in the Appendix. 

\subsubsection{}

Quasimodular forms occur in nature, for example, as
coefficients of the expansion of the odd 
genus $1$ theta-function
$$
\th(x) = (x^{1/2}-x^{-1/2}) \prod_{i=1}^\infty
 \frac{(1-q^i x)(1-q^i/x)}{(1-q^i)^2} 
$$
at the origin $x=1$. 
The techniques developed below give
a certain 
formula for \eqref{Z'} in terms of derivatives of
$\th(x)$ at $x=\pm 1$, from which the quasimodularity 
follows.

\subsubsection{}
The following discussion closely parallels the corresponding
discussion for the case of holomorphic differentials in \cite[\S{1.2}]{EO}.

Let $\sQ(\mu,\nu)$ denote the moduli space of pairs $(\Sigma, \phi)$, 
where $\phi$ is a quadratic differential on a curve $\Sigma$ with 
zeroes of multiplicities $\{\nu_i -2, 2 \mu_i - 2 \}$. Note that
we allow $\nu_i = 1$ hence our quadratic differentials can have simple
poles. For $(\Sigma, \phi) \in \sQ(\mu,\nu)$, let $\widetilde{\Sigma}$
denote the double cover of $\Sigma$ on which the differential 
$$
\omega=\sqrt{\phi}
$$
is well defined. The pair $(\widetilde{\Sigma}, \omega)$ belongs to
the corresponding space of holomorphic differentials
with zeroes of multiplicity
$$
\{\nu_i-1,\mu_i-1,\mu_i-1\} \,. 
$$
By construction, $\Sigma$ is the quotient of $\widetilde{\Sigma}$ by an
involution $\sigma$. Let $P$ denote the set of zeroes of $\omega$; it
is clearly stable under $\sigma$. Then $\sigma$ acts as an involution
on the relative homology group $H_1(\widetilde{\Sigma}, P, \Z)$. Let
$H^-$ denote the subspace of $H_1(\widetilde{\Sigma}, P, \Z)$ on which
$\sigma$ acts as multiplication by $-1$. Choose a basis 
$\{\gamma_1, \dots, \gamma_n\}$ for $H^-$, and consider the period map
$\Phi:  \sQ(\mu,\nu) \to \C^n$ defined by
\begin{displaymath}
\Phi(\Sigma,\phi) = \left( \int_{\gamma_1} \omega, \dots, \int_{\gamma_n}
  \omega \right).
\end{displaymath}
It is known \cite{MZ} that $\Phi(\Sigma,\phi)$ is a local
coordinate system on $\sQ(\mu,\nu)$ and, in particular, 
$n = \dim_{\C} H^- = \dim_{\C} \sQ(\mu,\nu)$. 

Pulling back the Lebesgue measure from $\C^n$ yields a well-defined
measure on $\sQ(\mu,\nu)$. However, this measure is infinite
since $\phi$ can be multiplied by any complex number. Thus we define
$\sQ_1(\mu,\nu)$ to be the subset satisfying 
\begin{displaymath}
\operatorname{Area}(\widetilde{\Sigma}) \equiv \frac{\sqrt{-1}}{2}
\int_{\widetilde{\Sigma}} \omega \wedge \overline{\omega} = 2
\end{displaymath}
As in the case of holomorphic differentials, the area function is a
quadratic form in the local coordinates on $\sQ(\mu,\nu)$, and
thus the image under $\Phi$ 
of $\sQ_1(\mu,\nu)$ can be identified with an open
subset of a hyperboloid in $\C^n$. 

Now let $E \subset \sQ_1(\mu,\nu)$ be a set lying in the domain
of a coordinate chart, and let $C \Phi(E) \subset \C^n$ 
denote the cone over $\Phi(E)$ with vertex $0$.  Then,
we can define a measure $\rho$ on $\sQ_1(\mu,\nu)$ via
\begin{displaymath}
\rho(E) = \operatorname{vol}(C \Phi(E)),
\end{displaymath}
where $\operatorname{vol}$ is the Lebesgue measure. 
The proof of \cite[Proposition 1.6]{EO} shows the following analogue:
\begin{displaymath}
\rho(\sQ_1(\mu,\nu)) = \lim_{D \to \infty} D^{-\dim_{\C} {\cal
    Q}(\mu,\nu)} \sum_{d=1}^{2D} \operatorname{Cov}^0_d(\mu,\nu),
\end{displaymath}
where $\operatorname{Cov}^0_d(\mu,\nu)$ is the number of inequivalent
degree $d$ connected covers $\cC \to \ffP$. Thus, the volume 
$\rho(\sQ_1(\mu,\nu))$ can be read off from the $q \to 1$
asymptotics of the connected generating function $Z^\circ(\mu,\nu; q)$.

Note that the moduli spaces $\sQ(\mu,\nu)$ may be disconnected. 
Ergodic theory applications require the knowledge of volumes of 
each connected component. Fortunately, connected components of 
$\sQ(\mu,\nu)$ have been classified by E.~Lanneau \cite{L} and these 
spaces turn out to be connected except for hyperelliptic 
components (whose volume can be computed separately) and 
finitely many sporadic cases.

\subsubsection{}

The modular transformation
$$
q=e^{2\pi i \tau} \mapsto e^{-2\pi i /\tau}
$$
relates $q=0$ and $q=1$ and thus gives an easy handle on
the $q\to 1$ asymptotics of \eqref{Z'}. This gives 
an asymptotic enumeration of pillowcase covers and, hence, 
computes the volume of the moduli spaces of quadratic 
differentials.

\subsubsection{}

In spirit, Theorem \ref{T1} is parallel to the results of 
\cite{BO,EOP,KaZ}, see also \cite{Dij,Dij2,Doug} for earlier 
results in the physics literature. 
The main novelty is the occurrence of 
quasimodular forms of higher level. One might speculate
whether similar lattice point enumeration in the space of
$N$th order differentials leads to level $N$ quasimodular
forms. Those spaces, however, do not admit an $SL_2(\R)$-action 
and a natural interpretation of their volumes is not known.

\subsubsection{}
The following enumerative problem is naturally a building 
block of the enumerative problem that we consider. 
Consider branched covers of the sphere ramified over 
3 points $0,1,\infty$ with profile 
$(\nu,2^{d-|\nu|/2})$, $(2^d)$, and $\mu$, respectively, 
where $\mu$ is an arbitrary partition of $2d$. 

The preimage of the segment $[0,1]$ on the sphere 
is a graph $\cG$ on a
Riemann surface (also known as a \emph{ribbon graph})
with a plenty of $2$-valent vertices (that can be ignored) 
and a few odd valent vertices (namely, with valencies $\nu_i$). 
The complement of $\cG$ is a union of $\ell(\mu)$ 
disks (known as \emph{cells}) with perimeters $2\mu_i$ 
in the natural metric on $\cG$. The asymptotic
enumeration of such combinatorial 
objects is, almost by definition, given by 
integrals of $\psi$-classes 
against Kontsevitch's combinatorial cycles 
in $\overline{\cM}_{g,\ell(\mu)}$, see \cite{Kon}. 
There is a useful expression for these integrals
in terms of Schur $Q$-functions obtained in 
\cite{DIZ,Joz}.  In fact, our original approach 
to the results presented in this paper 
was based on these ideas. 

While the proof that we 
give here is more direct, it is still interesting to 
investigate the connection with combinatorial classes
further, especially since a natural geometric interpretation 
of combinatorial classes is still missing. Perhaps the 
Gromov-Witten theory of the orbifold $\ffP$ is the natural 
place to look for it. This will be 
further discussed in \cite{OP2}.

\section{Character sums}

\subsection{Characters of near-involutions}

\subsubsection{}

There is a classical way to enumerate branched
coverings in terms of irreducible characters which 
is reviewed, for example, in \cite{J} or in \cite{OP}. Specialized
to our case, it gives
\begin{equation}
  \label{cov_char}
   Z(\mu,\nu;q) = \sum_{\lambda} q^{|\lambda|/2} \, 
\left(\frac{\dim \lambda}{|\lambda|!}\right)^2 \, 
\fb_{\nu,2,2,\dots}(\lambda) \, \fb_{2,2,\dots}(\lambda)^3 \, 
\prod_i \fb_{\mu_i}(\lambda) 
\end{equation}
where summation is over all partitions, $\dim\lambda$ is the dimension 
of the corresponding representation of the symmetric group, and 
$\fb_\eta(\lambda)$ is the \emph{central character} of an element 
with cycle type $\eta$ in representation $\lambda$. Recall that the 
sum of all permutations with cycle type $\eta$ acts as a scalar
operator in any representation $\lambda$ and, by definition, this 
number is $\fb_\eta(\lambda)$. In \eqref{cov_char}, 
as usual, 
we abbreviate $\fb_{k,1,1,\dots}$ to $\fb_k$.

\subsubsection{}
A lot is known about the characters of the symmetric 
group $S(2d)$ in the 
situation when the representation is arbitrary but the 
support of the permutation is bounded by some number 
independent of $d$. In particular, explicit formulas
exist for the functions $\fb_k$. 
 
Understanding the function $\fb_{\nu,2,2,\dots}$ is 
the key to evaluation of \eqref{cov_char}. 
That is, we must study characters of 
permutations that are a product of a permutation with 
finite support and a fixed-point-free involution. 
We call such permutations \emph{near-involutions}.

\subsubsection{}

By a result of Kerov and Olshanski \cite{KO},
the functions $\fb_k$ belong to the 
algebra $\Lambda^*$ generated by 
\begin{equation}
  \label{p_k}
    \pb_k(\lambda) = (1-2^{-k}) \zeta(-k) + \sum_{i}
\left[(\lambda_i-i+\tfrac12)^k -(-i+\tfrac12)^k\right]\,,
\end{equation}
moreover $\fb_k$ has weight $k+1$ in the weight filtration 
on $\Lambda^*$ defined by setting 
$$
\wt \pb_k = k+1 \,.
$$
The functions $\pb_k$ are central characters of
certain distinguished elements in the group 
algebra of symmetric group known as 
\emph{completed cycles}. 
See \cite{OP} for the discussion of the relation between $\pb_k$ and
$\fb_k$ from the viewpoint of Gromov-Witten theory.

\subsubsection{}

Our next goal is to generalize the results of \cite{KO}
to characters of near-involutions. This will require
enlarging the algebra of functions. In addition to 
the polynomials $\pb_k$, we will need quasi-polynomial
functions $\pbb_k$
defined in \eqref{def_pbk} below.

It is convenient to work with 
the following generating function 
$$
\eb(\lambda,z)\overset{\textup{def}}= \sum_i e^{z(\lambda_i-i+\frac12)} = 
\frac1z+ \sum_k  \pb_k(\lambda) \, \frac{z^k}{k!} \,.
$$
By definition, set 
\begin{align}\label{def_pbk}
\pbb_k(\lambda) &= i \, k ! [z^k] \,  \eb(\lambda,z+\pi i) \\
                  &= \sum_{i}
\left[(-1)^{\lambda_i-i+1}(\lambda_i-i+\tfrac12)^k -
(-1)^{-i+1}(-i+\tfrac12)^k\right] + \textup{const}\,,
\notag
\end{align}
where the constant terms are determined by the expansion 
$$
\sum_k \frac{z^k}{k!} \, \pbb_k(\emptyset)  = \frac{1}{e^{z/2}+e^{-z/2}} \,.  
$$
Up to powers of $2$, they are Euler numbers.

\subsubsection{}
Define
$$
\Lb = \Q\left[\pb_k,\pbb_k\right]_{k\ge 1} \,. 
$$
Setting 
$$
\wt \pbb_k  = k 
$$
gives the algebra $\Lb$ the weight grading. 
Note that 
if $f$ is homogeneous then 
\begin{equation}
  \label{parit}
   f(\lambda')=(-1)^{\wt f} \, f(\lambda)\,, 
\end{equation}
where $\lambda'$ denotes the conjugate partition.

\subsubsection{}\label{2-quo}
In the definition of $\Lb$, we excluded
the function 
$$
\pbb_0(\lambda) = \frac12 + \sum_{i}
\left[(-1)^{\lambda_i-i+1}-(-1)^{-i+1}\right]\,,
$$
which 
measures the difference between number of even and odd
numbers among $\{\lambda_i-i+1\}$, also known as the 
\emph{2-charge} of a partition $\lambda$. 

Every partition
$\lambda$ uniquely defines two partitions $\alpha$ and 
$\beta$, known as its \emph{2-quotients}, such that
$$
\left\{\lambda_i-i+\tfrac12\right\} = 
\left\{2(\alpha_i-i+\tfrac12) + \pbb_0(\lambda) \right\} \, 
\sqcup \, 
\left\{2(\beta_i-i+\tfrac12) - \pbb_0(\lambda) \right\} \,.
$$
A partition $\lambda$ will be called \emph{balanced} if 
$\pbb_0(\lambda) = \frac12$. 

Several constructions related to $2$-quotients will play an 
important role in this paper. A modern review of these ideas 
can be found, for example, in \cite{FL}. In 
particular, it is known that the character $\chi^\lambda_{2,2,\dots}$
of a fixed-point free involution in the representation 
$\lambda$ vanishes unless $\lambda$ is balanced, in 
which case 
\begin{equation}
  \label{chi222}
    \left|
\chi^\lambda_{2,2,\dots} \right|= \binom{|\lambda|/2}{|\alpha|,|\beta|} \, 
\dim \alpha \, \dim \beta \,.
\end{equation}
It follows that only balanced partitions contribute to 
the sum \eqref{cov_char}.

\subsubsection{}
For a balanced partition $\lambda$,  define  
\begin{equation}
  \label{def_gnu}
    \gb_\nu(\lambda)  = \frac{\fb_{(\nu,2,2,\dots)}(\lambda)}
{\fb_{(2,2,\dots)}(\lambda)}\,. 
\end{equation}
We will prove that this function lies in $\Lb$ in the following sense: 

\begin{Theorem}\label{T2}
The ratio \eqref{def_gnu} is the restriction of a
unique function $\gb_\nu\in\Lb$ of weight $|\nu|/2$
to the set of balanced 
partitions. 
\end{Theorem}

Several examples of the polynomials $\gb_\nu$ can be
found in the Appendix.

\subsubsection{}

In view of Theorem \ref{T2}, it is natural 
to introduce the following \emph{pillowcase 
weight}
$$
\bW(\lambda) = \left(\frac{\dim \lambda}{|\lambda|!}\right)^2\, 
 \fb_{2,2,\dots}(\lambda)^4 \,.
$$
Theorem \ref{T1} follows from \eqref{cov_char}, 
Theorem \ref{T2}, and 
the following result: 

\begin{Theorem}\label{T3} For any $F\in \Lb$, the following average
\begin{equation}
  \label{aver}
\lang F \rang_\bW = 
    \frac1{Z(\emptyset,\emptyset;q)}
\sum_{\lambda} q^{|\lambda|} \, \bW(\lambda) \, F(\lambda)
\end{equation}
is a polynomial in $E_2(q^2)$, $E_2(q^4)$, and $E_4(q^4)$
of weight $\wt F$. 
\end{Theorem}

Note that if $F$ is homogeneous of odd weight then 
$\lang F \rang_\bW=0$. This can be seen directly 
from \eqref{parit}. Also note that \eqref{aver}
will \emph{not} in general be of pure weight even 
if $F$ is a monomial in the 
generators $\pb_k$ and $\pbb_k$. This contrast 
with \cite{BO,EOP} hints to the existence of a 
better set of generators of the algebra $\Lb$. 
Probably such generators are related to 
descendents of orbifold points in the Gromov-Witten 
theory of $\ffP$. 

\subsubsection{}

It will be convenient to work with the following 
generating functions for the sums \eqref{aver}
\begin{equation}
  \label{npt}
  F(x_1,\dots,x_n) = \lang \prod 
\eb(\lambda, \ln x_i) \rang_\bW\,. 
\end{equation}
The function \eqref{npt} will be called the 
\emph{$n$-point function}. 

\subsection{Proof of Theorem \ref{T2}}

\subsubsection{}

In the proof of theorems \ref{T2} and \ref{T3} it will be 
very convenient to use fermionic Fock space formalism.
This formalism is standard and \cite{Kac,MJD} can be recommended
as a reference. A quick review of these techniques can 
be found, for example, in Section 2 of \cite{OP}. We follow the 
notation of \cite{OP}.

\subsubsection{}
By definition, the space $\LV$ is spanned by the following 
infinite wedge products
\begin{equation}
  \label{v_lam}
    v_\lambda = \underline{\lambda_1 - \tfrac12} \wedge
 \underline{\lambda_2 - \tfrac32} \wedge 
\underline{\lambda_3 - \tfrac52}\wedge \dots\,,
\end{equation}
where $\underline{k}$, $k\in \Z+\frac12$ is a basis of 
the underlying space $V$ and $\lambda$ is a partition. 
The subscript $0$ in $\LV$ refers to the charge zero 
condition: the $i$th factor in \eqref{v_lam} is 
$\underline{-i+\tfrac12}$ for all sufficiently large 
$i$. 

There is a natural projective representation of the 
Lie algebra $\gl(V)$ on $\LV$. For us, the following 
elements of $\gl(V)$ will be especially important:
\begin{equation}
  \label{Ef}
    \cE_k\left[f(x)\right] \, \underline{i} = 
f(i-\tfrac k2) \, \underline{i-k}\,,
\end{equation}
where $f$ is a function on the real line. To define
the action of $\cE_0\left[f(x)\right]$ on $\LV$ one 
needs to regularize the infinite sum $\sum_{i<0} f(\tfrac12-i)$. 
This regularization is the source of the central extension in 
the $\gl(V)$ action. When $f$ is an exponential as in 
$$
\cE_k(z) = \cE_k\left[e^{zx}\right]\,,
$$
this infinite sum is a geometric series and thus 
has a natural 
regularization. By differentiation, this leads to the 
$\zeta$-regularization for operators $\cE_k[f]$ with 
a polynomial function $f$.

\subsubsection{}

Other very useful operators are 
$$
\alpha_k = \cE_k[1] \,, \quad k\ne 0 \,. 
$$
The operator $H$ defined by 
$$
H \, v_\lambda  = |\lambda| \, v_\lambda
$$
is known as the energy operator. It differs only by a
constant from the operator $\cE_0[x]$. 
The operator $H$ defines a natural grading on 
$\LV$ and $\gl(V)$. 

\subsubsection{}

A function $F(\lambda)$ on partitions of $n$ can be viewed
as a vector 
$$
\sum_{|\lambda|=n} F(\lambda) \, v_\lambda \in \LV
$$
of energy $n$. For example, the vectors
\begin{equation}
  \label{char_v}
    \left| \mu \right\rangle \overset{\textup{def}}=
    \frac1{\zz(\mu)} \prod \alpha_{-\mu_i} \vac = 
\frac1{\zz(\mu)} 
\sum_{\lambda} \chi^\lambda_\mu \, v_\lambda \,,
\end{equation}
where
$$
\zz(\mu)=|\Aut \mu| \prod \mu_i  \,,
$$
correspond to irreducible characters normalized
by the order of the centralizer.

\subsubsection{}

The operator $\cE_0(z)$ is the generating function
$$
\cE_0(z) = \cE_0[e^{zx}]= \frac1z + \sum_{k} \frac{z^k}{k!} \, \cP_k \,, 
$$
for the operators $\cP_k$ acting by 
$$
\cP_k \, v_\lambda  = \pb_k(\lambda) \, v_\lambda  \,.
$$
In parallel to \eqref{def_pbk}, we define operators $\cPb_k$ by
$$
i \, \cE_0(z+\pi i) = \sum_{k} \frac{z^k}{k!} \, \cPb_k \,.
$$
Translated into the operator language, the
the statement of Theorem \ref{T2} is the following:
the orthogonal projection of $\ve{\nu,2^{d-|\nu|/2}}$ 
onto the subspace spanned by the $v_\lambda$ with $\lambda$ balanced
is a linear combination of vectors 
\begin{equation}
  \label{PP2}
    \prod \cP_{\mu_i}
\prod \cPb_{\bar\mu_i} \, \ve{2^{d}} 
\end{equation}
with 
$$
\wt \mu + |\bar\mu|  \le  |\nu|/2 
$$
and coefficients independent of $d$.

\subsubsection{}

Let us call the span of $v_\lambda$ with $\lambda$ 
balanced the balanced subspace of $\LV$. 
A convenient orthogonal basis of it 
is provided by the following vectors
\begin{equation}
  \label{char_sign}
    \left| \rho;\bar\rho \right\rangle \overset{\textup{def}}=
    \frac1{\zz(\rho)\zz(\bar\rho)} \prod \alpha_{-\rho_i} 
\prod \bar\alpha_{-\bar\rho_i}\,\vac \,,\quad  \rho_i,\bar\rho_i\in 2\Z\,,
\end{equation}
where the operators $\bar\alpha_k$ are defined by
\begin{equation}
  \label{def_alb}
  \bar\alpha_k = i^{k+1} \cE_k(\pi i) = \sum_n (-1)^{n+\frac12} 
E_{n-k,n} + \frac{\delta_{k}}{2} \,, 
\end{equation}
the operators $E_{i,j}$ being the matrix units of $\mathfrak{gl}(V)$. 
{F}rom the commutation relations for the operators $\cE_k(z)$ we 
compute
\begin{align}
[\bar\alpha_k,\bar\alpha_m] &= \left[(-1)^k - (-1)^m\right] \alpha_{k+m} +
k (-1)^k \delta_{k+m} \,, \label{com1}
\\
[\alpha_k,\bar\alpha_m] &= \left[1 - (-1)^k\right] 
\left(\bar\alpha_{k+m} +
\frac{\delta_{k+m}}{2}\right) \,. \label{com2}
\end{align}
In particular, when both $k$ and $m$ are even, all these operators
commute apart from the central term in $[\bar\alpha_k,\bar\alpha_{-k}]$. 

The adjoint of $\bar\alpha_k$ is 
$$
\bar\alpha_k^* = (-1)^k \bar\alpha_{-k} \,,
$$
which gives the following inner products 
\begin{equation}
  \label{inner}
     \lang \rho;\bar\rho \big| \rho';\bar\rho' \rang =
\frac{\delta_{\rho,\rho'} \, \delta_{\bar\rho,\bar\rho'}}
{\zz(\rho) \, \zz(\bar\rho)} \,,
\end{equation}
provided all parts of all partitions in \eqref{inner} are even. 
In particular, the vectors \eqref{char_sign} are orthogonal. 
It is clear that they lie in the balanced subspace and their 
number equals the dimension of the space. Therefore, they form 
a basis.


\subsubsection{}

The projection of $\ve{\nu,2^{d-|\nu|/2}}$ onto the balanced
subspace is given in term of inner products of the form
$$
 \lang \nu,2^{d-|\nu|/2} \big| (\rho,2^{d-|\rho|/2-|\bar\rho|/2}); \bar\rho \rang
$$
where all parts of $\nu$ are odd, 
all parts of $\rho$ and $\bar\rho$ are even, and 
$\rho$ has no parts equal to $2$. {F}rom the commutation 
relations \eqref{com1} and \eqref{com2} we conclude 
that this inner product vanishes unless
$$
\rho =\emptyset \,.
$$ 
The nonvanishing inner products are the following 
\begin{equation}\label{inn2}
\lang \nu,2^{k} \big| 2^k; \bar\rho \rang  = 
\frac{2^{\ell(\nu)-\ell(\bar\rho)}}{2^k k! \zz(\nu) \zz(\bar\rho)}\,
\ttC(\nu,\bar\rho)\,,
\end{equation}
where the combinatorial coefficient $\ttC(\nu,\bar\rho)$ equals the 
number of ways to represent the parts of $\bar\rho$ as sums of parts
of $\nu$. For example
$$
\ttC((3,1,1,1),(4,2)) = 3 \,, \quad \ttC((3,1,1,1),(6)) = 1 \,. 
$$


\subsubsection{}

The matrix elements 
\begin{equation}
  \label{mPP}
    \lang  2^d \left|\, \prod \cP_{\mu_i}
\prod \cPb_{\bar\mu_i} \right| (\rho,2^{d-|\rho|/2-|\bar\rho|/2}); \bar\rho \rang 
\,, \quad \rho_i\ne 2 \,, 
\end{equation}
describe the decomposition of the vectors \eqref{PP2} in 
the basis \eqref{char_sign}. Since 
\begin{equation}
  \label{cP_1}
  \cP_1 \, \ve{2^d} = (d-\tfrac{1}{24}) \, \ve{2^d} \,,
\end{equation}
we can also assume that $\mu_i\ne 1$. 

We claim that \eqref{mPP}
vanishes unless 
\begin{equation}
  \label{bound_wt}
   \wt \mu + |\bar\mu |
\ge  
\wt {\rho}/2 + |\bar\rho|/{2} \,,
\end{equation}
where $\rho/2$ is the partition with parts $\rho_i/2$ 
(recall that all parts of $\rho$ are even).

\subsubsection{}

The usual way to evaluate a matrix element like \eqref{mPP}
is to use commutation relations 
to commute all lowering operators to the right until 
they reach the vacuum (which they annihilate) and, 
similarly, commute the raising operators to the left. 

We will exploit the following property of the 
operators $\cP_k$ and $\cPb_k$: their commutator
with enough operators of the form $\alpha_{-2\rho_i}$
and $\bar\alpha_{-2\bar\rho_i}$ vanishes. All 
such commutators have the form $\cE_k[f]$ with 
$f(x)=(\pm 1)^x p(x)$, where $p(x)$ is a polynomial. 
Commutation with $\alpha_{-2\rho_i}$ takes a 
finite difference of $p(x)$; commutation with 
$\bar\alpha_{-2\bar\rho_i}$ additionally flips
the sign of $\pm 1$.

Since a $(k+1)$-fold finite difference of a degree
$k$ polynomial vanishes, the commutator of 
of $\cP_k$ with more than $k+1$ operators of the 
form $\alpha_{-2\rho_i}$ or $\bar\alpha_{-2\bar\rho_i}$
vanishes. In fact, a $(k+1)$-fold commutator may 
be nonvanishing only because of the central extension
term. To pick up this central term, the total 
energy of all operators involved should be zero and 
the number of $\bar\alpha$'s should be even. 
Same reasoning applies to $\cPb_k$, but now the 
the number of $\bar\alpha$'s should be odd 
to produce a nontrivial $(k+1)$-fold commutator.

\subsubsection{}

Now look at one of the raising operators involved in \eqref{mPP}, 
say $\bar\alpha_{-\bar\rho_i}$.  
This operator commutes with $\alpha_2$ and 
its adjoint annihilates the vacuum, so only the terms involving 
the commutator of $\bar\alpha_{-\bar\rho_i}$ with one of the $\cP_{\mu_i}$
or $\cPb_{\bar\mu_i}$ give a nonzero contribution to \eqref{mPP}. 
The commutator $[\cP_{\mu_i},\bar\alpha_{-\bar\rho_i}]$ has energy $(-\rho_i)$ 
and so its adjoint again annihilates the vacuum. The same is true for
the commutation with $\cPb_{\bar\mu_i}$.
To bring these commutators back to zero energy, 
one needs to commute it $\bar\rho_i/2$ times with $\alpha_2$. 
Given the above bounds on how many commutators we can 
afford, this implies \eqref{bound_wt}.

\subsubsection{}

When the bound \eqref{bound_wt} is saturated, then a further condition 
$$
\ell(\rho)+\ell(\bar\rho) \ge \ell(\mu)+\ell(\bar\mu)
$$
is clearly necessary for nonvanishing of \eqref{mPP}. The unique
nonzero coefficient saturating both bounds corresponds to 
$$
\rho = 2\mu\,, \quad \bar\rho = 2\bar\mu \,.
$$
Moreover, when divided by the
norm squared 
of the vector $\ve{(\rho,2^{d-|\rho|/2-|\bar\rho|/2}); \bar\rho}$, 
this coefficient is independent of $d$.

\subsubsection{}

For general $\rho$ and $\bar\rho$, the similarly normalized
coefficient will be a polynomial in $d$ of degree
\begin{equation}
  \label{wt_differ}
     \frac12\left(\wt \mu + |\bar\mu| - 
\wt \rho/2 - |\bar\rho|/2\right) 
\end{equation}
because so many operators $\alpha_{-2}$ can commute with $\cP_{\mu_i}$'s or $\cPb_{\bar\mu_i}$'s
instead of commuting directly with $\alpha_2$'s. 

By induction on weight and length, 
we can express the basis vectors \eqref{char_sign}
in terms of \eqref{PP2} with $\mu_i\ne 1$ and coefficients being 
polynomial in $d$ of degree at most minus the difference \eqref{wt_differ}. 
By \eqref{cP_1}, to have $d$-dependent coefficients and $\mu_i\ne 1$ is the same as
 to allow $\mu_i=1$ and make coefficients independent of $d$. The bound of degree in 
$d$ ensures that this transition preserves weight. 
This concludes the proof of 
Theorem \ref{T2}.

\section{Proof of Theorem \ref{T3}}

\subsection{The pillowcase operator}\label{sPO}

\subsubsection{}

Consider the following operator 
\begin{equation}
  \label{oW}
    \oW = 
\exp\left(\sum_{n>0} \frac{\alpha_{-2n-1}}{2n+1}\right)
\exp\left(-\sum_{n>0}
 \frac{\alpha_{2n+1}}{2n+1}\right)\,.
\end{equation}
Because this operator is normally ordered, its
matrix elements $(\oW v,w)$ are well defined
for any vectors $v$ and $w$ of finite energy. 
The relevance of this operator for our 
purposes lies in the following 

\begin{Theorem}\label{T4} 
The diagonal matrix elements of $\oW$ are
the following:  
$$
\left(\oW \, v_\lambda,v_\lambda\right) = 
\begin{cases}
 \bW(\lambda)\,, & \textup{$\lambda$ is balanced}\\
0\,, & \textup{otherwise} \,.
\end{cases}
$$
\end{Theorem}

The proof of this theorem will occupy the rest of 
Section \ref{sPO}.

\subsubsection{}

Let $N$ be chosen so large that $\lambda_{2N+1}=0$. 
Then because the operator \eqref{oW} is a product of 
an upper unitriangular and lower unitriangular operator, 
the vectors $\underline{\lambda_i-i+\tfrac12}$ with 
$i>2N$ in 
$$
v_\lambda = \underline{\lambda_1-\tfrac12} \wedge
\underline{\lambda_2-\tfrac32} \wedge 
\underline{\lambda_3-\tfrac52} \wedge \dots 
$$
are inert bystanders for the evaluation of 
$\left(\oW \, v_\lambda,v_\lambda\right)$. The whole
computation is therefore a computation of a matrix
element of an operator in a finite exterior power of 
a finite dimensional vector space 
$V^{[N]}$ with basis 
$$
e_k = \underline{-2N+k+\tfrac12}\,, \quad
k =0 \dots, \lambda_1+2N-1\,.
$$
By definition, matrix elements of $\oW$ in exterior powers
of $V^{[N]}$ are determinants of the matrix elements of 
$\oW$ acting on the space $V^{[N]}$ itself. The latter 
matrix elements are determined in the following 

\begin{Proposition}\label{P1} We have 
  \begin{equation}
    \label{oWme}
     \frac{\left(\oW \, e_k,e_l\right)}{\bb(k)\bb(l)} = 
\begin{cases}
1\,, & k\equiv l\equiv 0\mod 2\\
0\,, & k\equiv l\equiv 1\mod 2\\
2/(k-l)\,, & \textup{otherwise} \,, 
\end{cases}
\end{equation}
where
$$
\bb(k)= \frac{k!}{2^k \, \lfloor k/2 \rfloor!^2}\,. 
$$
\end{Proposition}

\subsubsection{}

For the proof of Proposition \ref{P1}, 
form the following generating function 
$$
f(x,y)=\sum_{k,l} x^k y^l \, \left(\oW \, e_k,e_l\right) \,.
$$
{F}rom the equality
$$
\exp\left(\sum_{n>0} \frac{x^{2n+1}}{2n+1}\right) = 
\sqrt{\frac{1+x}{1-x}}
$$
and definitions, we compute 
$$
f(x,y) =
\frac{1}{1-xy} \,  
\sqrt{\frac{1+x}{1-x}} \, \sqrt{\frac{1-y}{1+y}} 
\,.
$$
The following factorization 
$$
\left(x\frac{\partial}{\partial x} -
y\frac{\partial}{\partial y}\right) f(x,y)  =
\frac{(x+y)(1+x)(1-y)}{(1-x^2)^{3/2}(1-y^2)^{3/2}}\,,
$$
by elementary binomial coefficient manipulations
proves \eqref{oWme} for $k\ne l$. To compute the 
diagonal matrix elements observe that the 
above differential equation uniquely determines
$f(x,y)$ from its values on the diagonal $x=y$. 
On the diagonal, the skew-symmetric terms in
\eqref{oWme} cancel out and evaluation is 
immediate.

\subsubsection{}

We now proceed to the computation of the 
matrix element $\left(\oW \, v_\lambda,v_\lambda\right)$. 
We have the following 

\begin{Proposition} We have  
$$
\left(\oW \, v_\lambda,v_\lambda\right)  = 
\left(2^{N} 
\prod_{i=1}^{2N} \bb(\lambda_i-i+2N) \prod_{i < j \le 2N} 
(\lambda_i-\lambda_j+j-i)^{(-1)^{\lambda_i-\lambda_j+j-i}} \right)^2 \,,
$$
provided $\lambda$ is balanced and $\left(\oW \, v_\lambda,v_\lambda\right)=0$
otherwise. 
\end{Proposition}

The proof of this proposition is the following. Observe 
that by Proposition \ref{P1} the matrix element 
$\left(\oW \, v_\lambda,v_\lambda\right)$ is a determinant 
of a $2N\times 2N$ block matrix in which the odd-odd 
block is identically zero, the even-even block is 
a rank $1$ matrix with all elements equal to $1$ and 
the off-diagonal blocks 
have the form $\left(\frac2{x_i-y_j}\right)$, where 
$\{x_i\}$ and $\{y_i\}$ are the odd and even subsets
of $\{\lambda_i-i+2N\}$. Since the odd-odd block 
is identically zero, its size has to be $\le N$ for 
the determinant to be nonvanishing. Similarly, if the 
 size of the even-even block is larger than $N$ then 
the determinant is easily seen to vanish. It follows
that both blocks have size $N$, which precisely means 
that the partition $\lambda$ is balanced. It 
remains to use Cauchy determinant
$$
\det \left(\frac1{x_i+y_j}\right) = 
\frac{\prod_{i<j} (x_i-x_j) \, \prod_{i<j} (y_i-y_j)}
{\prod (x_i+y_j)}
$$ 
to finish the proof.

\subsubsection{}
Note that decomposition of $\{\lambda_i-i+2N\}$ into 
the even and odd subsets is  
the same as the $2$-quotient construction from 
Section \ref{2-quo}. 
Theorem \ref{T4} follows from formula \eqref{chi222}
and the following classical formula
$$
\frac{\dim \lambda}{|\lambda|!} = 
\frac{\prod_{i<j\le N} (\lambda_i-\lambda_j+j-i)}
{\prod (\lambda_i+N-i)!} \,,
$$
where $N$ is any number such that $\lambda_{N+1}=0$.

\subsubsection{}
It would be interesting to find an interpretation of
the operator $\oW$ in conformal field theory. Note that
$$
\exp\left(\sum_{n>0} \frac{z^{-2n-1}}{2n+1}\right)
\exp\left(-\sum_{n>0}
 \frac{z^{2n+1}}{2n+1}\right) = 
\sqrt{\frac{1+z^{-1}}{1-z^{-1}}} \, \sqrt{\frac{1-z}{1+z}}
$$
is the Wiener-Hopf factorization of the function taking the 
value $\mp i$ on the upper/lower half-plane.

\subsection{Formula for the $n$-point function }

\subsubsection{}

Theorem \ref{T4} yields the following operator 
formula for the $n$-point function \eqref{npt}
\begin{equation}
  \label{npt_trace}
  F(x_1,\dots,x_n) =   
\frac1{Z(\emptyset,\emptyset;q)} \tr 
q^H \prod \cE_0(\ln x_i) \, \oW \,,
\end{equation}
where the trace is taken in the charge zero 
subspace of the infinite wedge and $H$ is the 
energy operator
$$
H \, v_\lambda = |\lambda| \, v_\lambda \,.
$$
We have the following expression for the operator
$\cE_0$ in terms of the fermionic currents: 
$$
\cE_0(\ln x) = [y^0] \, \psi(xy) \, \psi^*(y)\,,
$$
where $[y^0]$ denotes the constant coefficient
in the Laurent series expansion in the variable $y$. 
Therefore, 
\begin{multline}\label{tr_psi}
 F(x_1,\dots,x_n) = \frac1{Z(\emptyset,\emptyset;q)} \times \\
[y_1^0 \dots y_n^0] \tr q^H \psi(x_1y_1) \, \psi^*(y_1) 
\dots \psi(x_n y_n) \, \psi^*(y_n) \, \oW\,. 
\end{multline}

\subsubsection{}
By the main result of \cite{FL}, we have 
\begin{equation}
  \label{bWle}
     \bW(\lambda) \le 1 
\end{equation}
for any partition $\lambda$. In other words,  all diagonal 
matrix elements of $\oW$ are bounded by $1$. For the 
off-diagonal elements, we prove the following cruder 
bound

\begin{Proposition} Let $M=\max\{|\lambda|,|\mu|\}$. Then
 \begin{equation}
    \label{bound}
      \left(\oW \, v_\lambda,v_\mu\right) \le 
\exp\left(\frac12 \sum_{n=0}^{\lfloor \frac{M-1}{2} \rfloor} \frac{1}{2n+1} 
\right) \sim \textup{const} \, M^{1/4} \,.
  \end{equation}
\end{Proposition}

To see this note that 
$$
 \left(\oW \, v_\lambda,v_\mu\right) =  
\left(\oW^{[M]} \, v_\lambda,v_\mu\right)
$$
where $\oW^{[M]}$ is the following truncated operator
$$
\exp\left(\sum_{2n+1 \le M} \frac{\alpha_{-2n-1}}{2n+1}\right)
\exp\left(-\sum_{2n+1 \le M}
 \frac{\alpha_{2n+1}}{2n+1}\right)\,.
$$
We claim that the operator $\oW^{[M]}$ is a multiple of 
a unitary operator. Indeed, 
$$
\left({\oW^{[M]}}^*\right)^{-1} = 
\exp\left(- \sum_{n=0}^{\lfloor \frac{M-1}{2} \rfloor} \frac{1}{2n+1} 
\right) \, \oW^{[M]}\,,
$$
whence the result. 

In fact, we will only use that \eqref{bound}
is bounded by a polynomial in the sizes of the partitions. 

\subsubsection{}

By normally ordering all 
fermionic operators in \eqref{tr_psi} and using the 
estimate  \eqref{bound} one sees that the 
trace converges if 
\begin{equation}\label{dom} 
  |y_n/q| > |x_1 y_1| > |y_1| > \dots > |x_n y_n| > |y_n| > 1 \,.
\end{equation}

\subsubsection{}
The proof of the following identity is given 
in \cite{Kac}, theorem 14.10: 
\begin{multline}
  \label{Kac}
  \psi(xy) \, \psi^*(y) = \frac1{x^{1/2}-x^{-1/2}} \times \\
\exp\left(\sum_n \frac{(xy)^n - y^n}{n}\, \alpha_{-n}\right) \,
\exp\left(\sum_n \frac{y^{-n} - (xy)^{-n}}{n}\, \alpha_{n}\right) \,. 
\end{multline}
It allows to express the operator in \eqref{tr_psi}
in terms of bosonic operators $\alpha_n$. 

With respect to the action of the operators $\alpha_n$, 
the charge zero subspace of the infinite wedge space 
decomposes as the following infinite tensor product 
$$
\LV \cong \bigotimes_{n=1}^\infty
 \bigoplus_{k=0}^\infty \alpha_{-n}^k \,\vac\,,
$$
the distinguished vector in each factor being $\vac$. 
This gives a factorization of the trace in \eqref{tr_psi}. 
The trace in each tensor factor is computed 
as follows:
$$
\tr e^{A \alpha_{-n}} \, e^{B \alpha_{n}} 
\Big|_{\bigoplus_{k=0}^\infty \alpha_{-n}^k \,\vac} = 
\frac1{1-q^n} \exp\left(\frac{nAB \, q^n}{1-q^n}\right) \,.
$$
For example, this shows that 
$$
\tr q^H \, \oW = (q^2)_\infty^{-1/2} = Z(\emptyset,\emptyset;q)\,,
$$
where
$$
(a)_\infty = \prod_{n\ge 0} (1-a q^n) \,,
$$ 
and so the $0$-point function is $F()=1$, as expected. For
the $n$-point function this gives the following 

\begin{Theorem} We have
\begin{multline}
  \label{nptf}
    F(x_1,\dots,x_n) = \prod \frac1{\th(x_i)} \times \\
[y_1^0\dots y_n^0] \,\,
\prod_{i<j} \frac{\th(y_i/y_j) \, \th(x_i y_i /x_j y_j)}
{\th(x_i y_i/y_j) \, \th(y_i /x_j y_j)} \, 
\prod_i \sqrt{\frac{\th(-y_i) \, \th(x_i y_i)}
{\th(y_i) \, \th(-x_i y_i)}}\,,
\end{multline}
where the series expansion is performed in the domain 
\eqref{dom}. 
\end{Theorem}

\subsection{Quasimodular forms}

\subsubsection{}

In the computation of \eqref{nptf}, we can 
assume that $1<|x_i|\ll |q^{-1}|$ for all $i$
and hence 
$$
|y_i| > |y_j| \prod |x_k|^{\pm 1} > |q y_i|\,, \quad i<j \,.
$$
The series expansion in \eqref{nptf} can then be 
performed using the following elementary

\begin{Lemma}\label{L1} We have
  \begin{equation}
    \label{int_theta}
      \frac1{2\pi i} \oint_{|y|=c} \frac{dy}{y} \, 
\prod_{i=1}^n \frac{\th(y/a_i)}{\th(y/b_i)} = 
\left(1-\prod\frac{a_i}{b_i}\right)^{-1} \, 
\sum_{i=1}^n \frac{\prod_j \th(b_i/a_j)}
{\prod_{j\ne i} \th(b_i/b_j)}\,,
\end{equation}
provided $c>|b_i|>|q| c$ for $i=1,\dots,n$. 
\end{Lemma}

This is obtained by computing the difference
of $\oint_{|y|=c}$ and $\oint_{|y|=|q|c}$ as a
sum of residues using 
$$
\th'(1)=1  \,. 
$$

\subsubsection{}

There are two obstacles to literary applying 
this lemma to the evaluation of \eqref{nptf}. 
The first are the square roots in \eqref{nptf}. 
However, we are ultimately interested in the
expansion of \eqref{nptf} about $x_i=\pm 1$. 
The expansion of the integrand 
about $x_i=\pm 1$ contains no square
roots, only theta function and its derivatives. 
Formulas for integrating derivatives can be 
obtained from \eqref{int_theta} by 
differentiating with respect to parameters.

\subsubsection{}\label{l'H}

The other issues is that at $x_i=1$ the 
integrand is an elliptic function of the 
corresponding $y_i$, and so the left hand side
of \eqref{int_theta} gives infinity times zero. 
This can be circumvented, for example, by 
replacing each factor of $x_i$ in the argument of 
each theta function an independent variable and
specializing them all back to $x_i$ only after
integration. By the L'Hospital rule, this 
will produce an additional differentiation 
any time we expand around $x_i=1$ for some $i$.

\subsubsection{}
In the end, we will get some rather complicated polynomial 
in theta functions and their derivatives
evaluated at $\pm 1$ divided by a power of
$\th(-1)$. This means that  we will get a combination of 
Eisenstein series arising from 
\begin{equation}
  \label{theta_1}
    \ln\frac{z}{\th(e^{z})} = 2 \sum_{k\ge 1} \frac{z^{2k}}{(2k)!} \, 
E_{2k}(q)\,,
\end{equation}
and
\begin{equation}
  \label{theta_-1}
    \ln\frac{\th(-e^{z})}{\th(-1)} = 
2 \sum_{k\ge 1} \frac{z^{2k}}{(2k)!} \, 
\left[E_{2k}(q)-2^{2k} E_{2k}(q^2)\right]\,, 
\end{equation}
together with the following product
\begin{equation}
  \label{Null}
  \th(-1) = 2i \left( \prod_n \frac{1+q^n}{1-q^n} \right)^2 =
\frac{\eta(q^2)^2}{\eta(q)^4} \,.
\end{equation}
Note that \eqref{Null} has weight $-1$. 

\subsubsection{}

Without knowing the precise form of the answer, 
one can still make some qualitative observations
about it. 

Suppose we are interested in the coefficient 
of $z_1^{k_1}\dots z_n^{k_n}$ in the expansion of 
$$
F(e^{z_1},\dots,e^{z_r},-e^{z_{r+1}},\dots,-e^{z_n}) 
$$
in powers of $z_i$. We claim that the weight of
this coefficient is at most $\sum k_i + r$. 
Indeed, we from \eqref{theta_1} and \eqref{theta_-1} we have 
$$
\wt \, \left(x \frac{d}{dx}\right)^k \, \th(x)
\Big|_{x=\pm 1} = k-1 \,.
$$
This gives the following count for the weight
$$
n - n + \sum k_i + r\,,
$$
where the first $n$ is added because of 
the prefactor in \eqref{nptf}, 
the second $n$ is subtracted due to integration in $y_i$ 
(which, by Lemma \ref{L1} changes the balance of
$\theta$-factors by $1$), $\sum k_i$ is the 
number of times we need to differentiate the integrand,
and, finally, $r$ additional differentiations are needed 
for reasons explained in Section \ref{l'H}.

\subsubsection{}

We further claim that \eqref{nptf} is, in fact, a 
polynomial in the coefficients of \eqref{theta_1}, 
\eqref{theta_-1} and 
\begin{equation}
  \label{Eis4}
  \frac1{\th(-1)^2} = - \frac14 
\frac{\eta(q)^8}{\eta(q^2)^4} = 2 E_2(q) -12  E_2(q^2) + 16 E_2(q^4)\,. 
\end{equation}
First observe only even powers of \eqref{Null} appear in 
the answer. This is because the formula \eqref{nptf} has
a balance of minus signs in the arguments of theta 
functions in the numerator and denominator. 
Every time we specialize $y_i$ to one of 
the poles in \eqref{int_theta}, 
the balance of minus signs changes by an even number. 

Inverse powers of \eqref{Eis4} cannot appear in the 
answer because they grow exponentially as $q\to 1$
and there are no other exponentially large terms
to cancel this growth out. The averages \eqref{aver} may grow only 
polynomially as $q\to 1$ because of the bound \eqref{bWle}. 

\subsubsection{}\label{squasi}

Recall from \cite{KZ} that a \emph{quasimodular
form} for a congruence subgroup $\Gamma\subset SL_2(\Z)$ 
is, by definition, the holomorphic part of an 
almost holomorphic modular form for $\Gamma$. 
A function of $|q|<1$ is called almost holomorphic
if it is a polynomial in $(\ln |q|)^{-1}$ with 
coefficients in holomorphic functions of $q$. 
Quasimodular forms for $\Gamma$ forms a graded
algebra denoted by $\QM(\Gamma)$. By a theorem of
Kaneko and Don Zagier \cite{KZ}, 
$$
\QM(\Gamma) = \Q[E_2] \otimes \cM(\Gamma) \,.
$$
In particular, 
\begin{equation}
  \label{3Eis}
     E_2(q), E_2(q^2), E_2(q^4) \in \QM(\Gamma_0(4))
\end{equation}
where 
$$
\Gamma_0(4)= \left\{
  \begin{pmatrix}
    a & b \\
  c & d
  \end{pmatrix}\,\Big|\,\, c\equiv 0 \bmod 4 
\right\} \subset SL_2(\Z) \,.
$$
Hence all averages \eqref{aver} lie in $\QM(\Gamma_0(4))$. 

In fact, the series \eqref{3Eis} generate
the subalgebra $\QM_{2*}(\Gamma_0(4))$ of 
even weight quasimodular forms. This is because
$\cM_{2*}(\Gamma_0(4))$ is freely generated by 
two generators of weight two, for example, by
 $E^\textup{odd}_2(q)$
and $E^\textup{odd}_2(q^2)$, where 
$$
E^\textup{odd}_2(q) =
E_2(q)-2E_2(q^2) = \frac1{24} + \sum_{n=1}^\infty \left(
\sum_{d|n, \,\textup{$d$ odd}} d\right) q^n \,. 
$$

\subsubsection{}
Note that because $\bW(\lambda)=0$ for any partition 
$\lambda$ of odd size, the series \eqref{aver} is in 
fact a series in $q^2$. It follows that it is 
quasimodular with respect to a bigger group, namely
$$
\begin{pmatrix}
   1 & 0 \\
  0 & 2
  \end{pmatrix} \, 
\Gamma_0(2) \, 
\begin{pmatrix}
    1 & 0 \\
  0 & 2
  \end{pmatrix}^{-1} \supset \Gamma_0(4) \,.
$$
In other words, \eqref{aver} is, in fact, obtained
by substituting $q\mapsto q^2$ into an element of 
$\QM(\Gamma_0(2))$. We have 
$$
\cM(\Gamma_0(2))= \Q[E^\textup{odd}_2(q),E_4(q^2)]
$$
and hence 
$$
\QM(\Gamma_0(2)) = \Q[E_2(q),E_2(q^2),E_4(q^2)] \,.
$$
This concludes the proof of Theorem \ref{T3}.

\appendix

\section{Examples}

In this appendix, we list some simple examples of the 
quasimodular forms $Z'(\mu,\nu;q)$ appearing in 
Theorem \ref{T1} and polynomials  $\gb_\nu$ 
from Theorem \ref{T2}.

\subsection{Quasimodular forms $Z'(\mu,\nu;q)$}

$$
Z'((1,1),(2))=20\,{{ E_2(q^4)}}^{2}-20\,{ E_2(q^4)}\,{ E_2(q^2)}+4\,{{ E_2(q^2)}}^{2}-5/3\,{
 E_4(q^4)}\,.
$$

\begin{multline*}
Z'((3,1),(3)) = 
-{\frac {2112}{5}}\,{{ E_2(q^4)}}^{3}+{\frac {3888}{5}}\,{{ E_2(q^4)}}^{2}
{ E_2(q^2)}-{\frac {2304}{5}}\,{ E_2(q^4)}\,{{ E_2(q^2)}}^{2}\\
+
{\frac {384}{5}}
\,{{ E_2(q^2)}}^{3}+48\,{ E_4(q^4)}\,{ E_2(q^4)}-36\,{ E_4(q^4)}\,{ E_2(q^2)}
\,.
\end{multline*}

\begin{multline*}
Z'((3,3),(2)) ={\frac {1056}{5}}\,{{ E_2(q^4)}}^{3}-{\frac {1044}{5}}\,{{ E_2(q^4)}}^{2}{
 E_2(q^2)}\\
+{\frac {252}{5}}\,{ E_2(q^4)}\,{{ E_2(q^2)}}^{2}
-{\frac {12}{5}}\,{
{ E_2(q^2)}}^{3}-24\,{ E_4(q^4)}\,{ E_2(q^4)}+3\,{ E_4(q^4)}\,{ E_2(q^2)}
\\ 
+15/2\,{{
 E_2(q^4)}}^{2}
-15/2\,{ E_2(q^4)}\,{ E_2(q^2)}+3/2\,{{ E_2(q^2)}}^{2}-5/8\,{ 
E_4(q^4)}\,.
\end{multline*}

\begin{multline*}
Z'((5,1),(2)) = {\frac {3520}{3}}\,{{ E_2(q^4)}}^{3}-1160\,{{ E_2(q^4)}}^{2}{ E_2(q^2)}
+280\,
{ E_2(q^4)}\,{{ E_2(q^2)}}^{2}\\
-{\frac {40}{3}}\,{{ E_2(q^2)}}^{3}-{\frac {400}
{3}}\,{ E_4(q^4)}\,{ E_2(q^4)}+{\frac {50}{3}}\,{ E_4(q^4)}\,{ E_2(q^2)}+{\frac {
125}{3}}\,{{ E_2(q^4)}}^{2}\\
-{\frac {125}{3}}\,{ E_2(q^4)}\,{ E_2(q^2)}+{
\frac {25}{3}}\,{{ E_2(q^2)}}^{2}-{\frac {125}{36}}\,{ E_4(q^4)}\,.
\end{multline*}

$$
Z'((1,1,1,1),\emptyset) = 1/4\,{ E_2(q^4)}+{\frac {1}{96}}\,.
$$

\begin{multline*}
Z'((3,3,1,1),\emptyset) = 
{\frac {9}{256}}-12\,{{ E_2(q^4)}}^{2}+{\frac {27}{2}}\,{ E_2(q^4)}\,{ 
E_2(q^2)}\\
-9/4\,{{ E_2(q^2)}}^{2}+5/4\,{ E_4(q^4)}+{\frac {9}{16}}\,{ E_2(q^4)}+3/8\,
{ E_2(q^2)}\,. 
\end{multline*}

\begin{multline*}
Z'((5,1,1,1),\emptyset) = 
{\frac {125}{1152}}-10\,{{ E_2(q^4)}}^{2}+15\,{ E_2(q^4)}\,{ E_2(q^2)}
\\
-5/2\,{
{ E_2(q^2)}}^{2}+{\frac {55}{24}}\,{ E_2(q^4)}+{\frac {5}{12}}\,{ E_2(q^2)}\,.
\end{multline*}

\begin{multline*}
Z'((3,3,3,3),\emptyset) =
-{\frac {24}{5}}\,{{ E_2(q^4)}}^{3}-{\frac {84}{5}}\,{{ E_2(q^4)}}^{2}{
 E_2(q^2)}+{\frac {423}{20}}\,{ E_2(q^4)}\,{{ E_2(q^2)}}^{2}
\\-{\frac {39}{10}}
\,{{ E_2(q^2)}}^{3}+{ E_4(q^4)}\,{ E_2(q^4)}+\frac74\,{ E_4(q^4)}\,{ E_2(q^2)}
-{\frac {
33}{4}}\,{{ E_2(q^4)}}^{2}
+{\frac {141}{16}}\,{ E_2(q^4)}\,{ E_2(q^2)}
\\
-{\frac {21}{32}}\,{{ E_2(q^2)}}^{2}+{\frac {25}{32}}\,{ E_4(q^4)}+{\frac {27}
{256}}\,{ E_2(q^4)}+{\frac {9}{32}}\,{ E_2(q^2)}+{\frac {27}{2048}}\,. 
\end{multline*}

\begin{multline*}
Z'((5,3,3,1),\emptyset) =
132\,{{ E_2(q^4)}}^{3}-708\,{{ E_2(q^4)}}^{2}{ E_2(q^2)}+639\,{ E_2(q^4)}\,{{
 E_2(q^2)}}^{2}\\
-114\,{{ E_2(q^2)}}^{3}-15\,{ E_4(q^4)}\,{ E_2(q^4)}+55\,{ E_4(q^4)}
\,{ E_2(q^2)}-310\,{{ E_2(q^4)}}^{2}\\
+{\frac {1365}{4}}\,{ E_2(q^4)}\,{ E_2(q^2)}
-{\frac {285}{8}}\,{{ E_2(q^2)}}^{2}+{\frac {175}{6}}\,{ E_4(q^4)}+{\frac {
615}{64}}\,{ E_2(q^4)}+{\frac {85}{8}}\,{ E_2(q^2)}+{\frac {375}{512}} \,. 
\end{multline*}

\subsection{Polynomials $\gb_\nu$}

\begin{equation*}
\gb_{1,1} = 1/2\,{\pbb}_{{1}}\,.
\end{equation*}

\begin{equation*}
\gb_{3,1} = 1/6\,{{\pbb}_{{1}}}^{2}+1/6\,{\pbb}_{{2}}-1/2\,\pb_{{1}}\,.
\end{equation*}

\begin{equation*}
\gb_{3,3} = 
-{\frac {1}{54}}\,{{\pbb}_{{1}}}^{3}+1/18\,{\pbb}_{{1}}{\pbb}_{{
2}}+{\frac {1}{54}}\,{\pbb}_{{3}}-1/4\,\pb_{{2}}+3/16\,{\pbb}_{{1}}
\,.
\end{equation*}

\begin{equation*}
\gb_{5,1} = 
1/30\,{{\pbb}_{{1}}}^{3}+1/10\,{\pbb}_{{1}}{\pbb}_{{2}}-1/2\,{
\pbb}_{{1}}\pb_{{1}}+1/15\,{\pbb}_{{3}}-1/2\,\pb_{{2}}+{\frac {25}{24}
}\,{\pbb}_{{1}}
\,.
\end{equation*}

\begin{multline*}
\gb_{5,3} = 
-{\frac {1}{360}}\,{{\pbb}_{{1}}}^{4}-{\frac {1}{60}}\,{{\pbb}_{{1
}}}^{2}{\pbb}_{{2}}-1/12\,{{\pbb}_{{1}}}^{2}\pb_{{1}}+{\frac {2}{45}
}\,{\pbb}_{{3}}{\pbb}_{{1}}+{\frac {25}{36}}\,{{\pbb}_{{1}}}^{2}
+1/40\,{{\pbb}_{{2}}}^{2}\\
-1/12\,{\pbb}_{{2}}\pb_{{1}}+5/8\,{\pb_{{1}}}
^{2}+{\frac {1}{60}}\,{\pbb}_{{4}}-1/2\,\pb_{{3}}+{\frac {25}{36}}\,{
\pbb}_{{2}}-{\frac {25}{12}}\,\pb_{{1}}
\,.
\end{multline*}

\begin{equation*}
\gb_{1,1,1,1} = 
-1/24\,{{\pbb}_{{1}}}^{2}+1/12\,{\pbb}_{{2}}+{\frac {1}{96}}
\,.
\end{equation*}

\begin{equation*}
\gb_{3,1,1,1} = 
{\frac {1}{108}}\,{{\pbb}_{{1}}}^{3}-1/36\,{\pbb}_{{1}}{\pbb}_{{
2}}-1/4\,{\pbb}_{{1}}\pb_{{1}}+{\frac {2}{27}}\,{\pbb}_{{3}}+3/8\,{
\pbb}_{{1}}
\,.
\end{equation*}

\begin{multline*}
\gb_{3,3,1,1} = 
{\frac {1}{216}}\,{{\pbb}_{{1}}}^{4}-1/12\,{{\pbb}_{{1}}}^{2}\pb_{{1
}}+{\frac {1}{108}}\,{\pbb}_{{3}}{\pbb}_{{1}}-1/8\,\pb_{{2}}{\pbb}
_{{1}}+{\frac {9}{32}}\,{{\pbb}_{{1}}}^{2}-{\frac {1}{72}}\,{{\pbb
}_{{2}}}^{2}\\
-1/12\,{\pbb}_{{2}}\pb_{{1}}+1/8\,{\pb_{{1}}}^{2}+1/36\,{
\pbb}_{{4}}+{\frac {9}{16}}\,{\pbb}_{{2}}-3/4\,\pb_{{1}}+{\frac {9}{
256}}
\,.
\end{multline*}

\begin{multline*}
\gb_{3,3,3,1} = 
{\frac {1}{4860}}\,{{\pbb}_{{1}}}^{5}+{\frac {1}{486}}\,{{\pbb}_{{
1}}}^{3}{\pbb}_{{2}}+{\frac {1}{108}}\,{{\pbb}_{{1}}}^{3}\pb_{{1}}-{
\frac {5}{972}}\,{\pbb}_{{3}}{{\pbb}_{{1}}}^{2}-1/24\,\pb_{{2}}{{
\pbb}_{{1}}}^{2}-{\frac {1}{96}}\,{{\pbb}_{{1}}}^{3}\\
+{\frac {1}{
324}}\,{\pbb}_{{1}}{{\pbb}_{{2}}}^{2}-1/36\,{\pbb}_{{1}}{\pbb}
_{{2}}\pb_{{1}}+{\frac {1}{162}}\,{\pbb}_{{4}}{\pbb}_{{1}}-{\frac {5
}{972}}\,{\pbb}_{{3}}{\pbb}_{{2}}-{\frac {1}{108}}\,{\pbb}_{{3}}
\pb_{{1}}-1/24\,\pb_{{2}}{\pbb}_{{2}}\\
+1/8\,\pb_{{2}}\pb_{{1}}+{\frac {31}{96
}}\,{\pbb}_{{1}}{\pbb}_{{2}}-{\frac {19}{32}}\,{\pbb}_{{1}}\pb_{{1
}}+1/4\,{\pbb}_{{3}}-\pb_{{2}}+{\frac {2}{405}}\,{\pbb}_{{5}}+{
\frac {153}{128}}\,{\pbb}_{{1}}
\,.
\end{multline*}

\begin{multline*}
\gb_{3,3,3,3} = 
{\frac {1}{29160}}\,{{\pbb}_{{1}}}^{6}-{\frac {1}{2916}}\,{\pbb}_{
{3}}{{\pbb}_{{1}}}^{3}+{\frac {1}{216}}\,\pb_{{2}}{{\pbb}_{{1}}}^{3}
-{\frac {1}{432}}\,{{\pbb}_{{1}}}^{4}+{\frac {1}{1944}}\,{{\pbb}_{
{1}}}^{2}{{\pbb}_{{2}}}^{2}-{\frac {1}{972}}\,{\pbb}_{{4}}{{\pbb
}_{{1}}}^{2}\\
+{\frac {1}{972}}\,{\pbb}_{{3}}{\pbb}_{{1}}{\pbb}_{{
2}}-{\frac {1}{72}}\,\pb_{{2}}{\pbb}_{{1}}{\pbb}_{{2}}-{\frac {7}{
288}}\,{{\pbb}_{{1}}}^{2}{\pbb}_{{2}}-1/12\,{{\pbb}_{{1}}}^{2}\pb_
{{1}}-{\frac {1}{2916}}\,{{\pbb}_{{2}}}^{3}-{\frac {1}{1944}}\,{{
\pbb}_{{3}}}^{2}\\
-{\frac {1}{216}}\,{\pbb}_{{3}}\pb_{{2}}+{\frac {59}
{864}}\,{\pbb}_{{3}}{\pbb}_{{1}}+1/32\,{\pb_{{2}}}^{2}-{\frac {3}{64
}}\,\pb_{{2}}{\pbb}_{{1}}+{\frac {1}{1215}}\,{\pbb}_{{5}}{\pbb}_{{
1}}+{\frac {231}{512}}\,{{\pbb}_{{1}}}^{2}+1/32\,{{\pbb}_{{2}}}^{2
}\\
-1/12\,{\pbb}_{{2}}\pb_{{1}}+3/8\,{\pb_{{1}}}^{2}+{\frac {5}{144}}\,{
\pbb}_{{4}}-{\frac {5}{12}}\,\pb_{{3}}+{\frac {1}{2916}}\,{\pbb}_{{6
}}+{\frac {129}{256}}\,{\pbb}_{{2}}-{\frac {9}{8}}\,\pb_{{1}}+{\frac {
27}{2048}}
\,.
\end{multline*}


\begin{thebibliography}{99}



\bibitem{BO}
S.~Bloch and A.~Okounkov,
\emph{The Character of the Infinite Wedge Representation},
Adv.\ Math.\ \textbf{149} (2000), no.~1, 1--60. 





\bibitem{Dij}
R.~Dijkgraaf,
\emph{Mirror symmetry and elliptic curves},
The Moduli Space of Curves, R.~Dijkgraaf,
C.~Faber, G.~van~der~Geer (editors), 
Progress in Mathematics, \textbf{129},
Birkh\"auser, 1995.



\bibitem{Dij2}
R.~Dijkgraaf,
\emph{Chiral deformations of conformal field theories},
Chiral deformations of conformal field theories, 
Nuclear Phys, B 493 (1997), no.\ 3, 588--612.

\bibitem{DIZ}
Ph.~Di Francesco, C.~Itzykson, J.-B.~Zuber, 
\emph{Polynomial averages in the Kontsevich model},
Comm.\ Math.\ Phys.\ \textbf{151} (1993), no. 1, 193--219.


\bibitem{Doug}
M.~Douglas,
\emph{Conformal field theory techniques in large 
$N$ Yang-Mills theory},
Quantum field theory and string theory (Carg\`ese, 1993),  
119--135, NATO Adv.\ Sci.\ Inst.\ Ser.\ B Phys., 328, Plenum, New York, 1995.


\bibitem{EMZ}
A.~Eskin, H.~Masur, and A.~Zorich, 
\emph{Moduli spaces of abelian differentials: 
the principal boundary, counting problems, 
and the Siegel-Veech constants},  
Publ.\ Math.\ Inst.\ Hautes \'Etudes Sci.  
\textbf{97} (2003), 61--179. 

\bibitem{EO}
A.~Eskin, A.~Okounkov, 
\emph{Asymptotics of numbers of branched coverings 
of a torus and volumes of moduli spaces of holomorphic differentials},
Invent.\ Math.\ \textbf{145} (2001), no.\ 1, 59--103. 


\bibitem{EOP}
A.~Eskin, A.~Okounkov, and R.~Pandharipande, 
\emph{The theta characteristic of a branched covering},
math.AG/0312186. 

\bibitem{FL}
S.~Fomin and N.~Lulov, 
\emph{On the number of rim hook tableaux},
Zap.\ Nauchn.\ Sem.\ S.-Peterburg.\ Otdel.\ Mat.\ Inst.\ Steklov
 (POMI) \textbf{223} (1995),
219--226, 
translation in J.\ Math.\ Sci.\ \textbf{87} (1997), no.~6, 4118--4123. 



\bibitem{J}
G.~Jones, 
\emph{Characters and surfaces: a survey} 
The atlas of finite groups: ten years on (Birmingham, 1995), 90--118, 
London Math.\ Soc.\ Lecture Note Ser., \textbf{249}, 
Cambridge Univ. Press, Cambridge, 1998. 


\bibitem{Joz}
T.~J\'ozefiak, 
\emph{Symmetric functions in the Kontsevich-Witten intersection theory of the moduli space of curves},
Lett.\ Math.\ Phys.\ \textbf{33} (1995), no.~4, 347--351.

\bibitem{Kac}
V.~Kac,
\emph{Infinite dimensional Lie algebras},
Cambridge University Press.


\bibitem{KaZ}
M.~Kaneko, D.~Zagier,
\emph{A generalized Jacobi theta function and quasimodular forms},
The moduli space of curves (Texel Island, 1994), 165--172, 
Progr.\ Math., 129, Birkhäuser Boston, Boston, MA, 1995.


\bibitem{KO}
S.~Kerov and G.~Olshanski,
\emph{Polynomial functions on the set of Young diagrams},
C.~R.~Acad.\ Sci.\
Paris S\'er.~I Math., \textbf{319}, no.~2, 1994, 121--126.



\bibitem{Kon}  
M.\ Kontsevich, {\em Intersection theory on the moduli
space of curves and the matrix Airy function}, Comm.\ Math.\ Phys.
{\bf 147} (1992), 1-23.


\bibitem{KZ}
M.~Kontsevich, A.~Zorich
\emph{Connected components of the moduli spaces 
of abelian differentials with prescribed singularities}, 
math.GT/0201292. 



\bibitem{L}
E.~Lanneau,
\emph{Classification des composantes connexes des strates de l'espace 
des modules des differentielles quadratiques}, 
Ph.D.\ Universit\'e Rennes I 2003. 


\bibitem{MZ}
H.~Masur and  A.~Zorich, 
\emph{Multiple Saddle Connections on Flat Surfaces and Principal Boundary of the Moduli Spaces of Quadratic Differentials}, 
math.GT/0402197. 

\bibitem{Mir}
M.~Mirzakhani, 
see {\tt http://math.harvard.edu/$\sim$mirzak/}. 

\bibitem{MJD}
T.~Miwa, M.~Jimbo, E.~Date, 
\emph{Solitons. 
Differential equations, symmetries and infinite-dimensional algebras},
Cambridge Tracts in Mathematics, 135. 
Cambridge University Press, Cambridge, 2000. 


\bibitem{OP}
A.~Okounkov and R.~Pandharipande,
\emph{Gromov-Witten theory, Hurwitz theory, and completed cycles},
math.AG/0204305.

\bibitem{OP2}
A.~Okounkov and R.~Pandharipande,
in preparation. 


\end{thebibliography}
\end{document}